\documentclass[12pt,a4paper,oneside,reqno]{article}

\usepackage[utf8]{inputenc}
\usepackage[english]{babel}

\usepackage{amsmath}
\usepackage{amssymb}
\usepackage{amsfonts}
\usepackage{amsthm}
\usepackage{stmaryrd}
\usepackage{cancel}
\usepackage{eurosym}
\usepackage{mathrsfs}

\usepackage[mathcal]{euscript}

\usepackage{tikz}
\usetikzlibrary{arrows}
\usepackage{tkz-graph}

\usepackage{longtable}
\usepackage{array}
\usepackage{booktabs}
\usepackage{caption}
\usepackage{tasks}
\usepackage{enumerate}
\usepackage{makecell}

\usepackage{cite}
\usepackage{url}
\usepackage{hyperref}
\usepackage{cleveref}

\usepackage{geometry}
\usepackage{fancyhdr} 
\usepackage{xcolor}

\usepackage{pdflscape}
\usepackage{ulem}
\usepackage{verbatim}
\usepackage{comment}
\usepackage[symbol]{footmisc}
\usepackage{float}

\theoremstyle{plain}
\newtheorem{theorem}{Theorem}

\newtheorem{remark}[theorem]{Remark}

\newtheorem{definition}[theorem]{Definition}
\newtheorem{corollary}[theorem]{Corollary}

\usepackage{fouriernc}

\begin{document}

\noindent{\Large Classification of three-dimensional Leibniz dialgebras\\ with two-dimensional annihilator}

\begin{center}

{\bf Naurizbay Uzakbaev\footnote{V.I.Romanovskiy Institute of Mathematics, Academy of Sciences of Uzbekistan, and Karakalpak State University, Uzbekistan; \ \ nn.uzakbaev@gmail.com}}

\end{center}

\noindent {\bf Abstract:} {\it A Leibniz dialgebra is a vector space with two bilinear products $\dashv$ and $\vdash$ satisfying the dialgebra analogue of the right Leibniz identity. This notion is obtained by the same standard procedure that turns the associative identity into the axioms of Loday's associative dialgebras. In particular, $(\mathcal{DL},\dashv)$ is a right Leibniz algebra, so a Leibniz dialgebra is a Leibniz algebra together with a second, compatible product. This paper recalls the definition of Leibniz dialgebras and continue the algebraic classification of three-dimensional complex Leibniz dialgebras. This paper describes up to isomorphism all Leibniz dialgebras $(\mathcal{DL},\dashv,\vdash)$ for which the Leibniz algebra $(\mathcal{DL},\dashv)$ is isomorphic to one of the algebras $\mathcal{L}_{8}\left(\alpha\right),$ $\mathcal{L}_9,$ $\mathcal{L}_{10},$ $\mathcal{L}_{11}$ of the classification of three-dimensional complex Leibniz algebras. The result consists $18$ families of pairwise non-isomorphic Leibniz dialgebras, nine of which depend on a complex parameter.}

\bigskip


\bigskip


\bigskip\bigskip

\tableofcontents

\newpage

\section{Introduction}

Leibniz algebras were introduced by Bloh~\cite{Bloh} under the name of $D$-algebras and were rediscovered by Loday~\cite{Loday93} as a non-antisymmetric analogue of Lie algebras: a (right) Leibniz algebra is an algebra in which every right multiplication is a derivation. Every Lie algebra is a Leibniz algebra, and a Leibniz algebra in which $\left[x, x\right]=0$ for all $x$ is a Lie algebra. Since then Leibniz algebras have been studied intensively from the algebraic, homological and geometric points of view. In particular, Loday and Pirashvili~\cite{LodayPirashvili} constructed the universal enveloping algebra and developed the (co)homology theory of Leibniz algebras; a systematic account of the structure theory and classification of Leibniz algebras is given in~\cite{AOR}. Among recent results, we mention the geometric classification of symmetric Leibniz algebras obtained in~\cite{FKS} and the description of solvable Leibniz superalgebras whose nilradical has the characteristic sequence $(n-1,1\,|\,m)$ and nilindex $n+m$ given in~\cite{KhM}.

Looking for an algebraic structure which plays for Leibniz algebras the role that associative algebras play for Lie algebras, Loday introduced (associative) dialgebras, also called diassociative algebras~\cite{Loday95,Loday01}. A dialgebra is a vector space with two associative products $\dashv$ and $\vdash$ satisfying three additional compatibility identities (see Definition~\ref{def:diassociative}), and the bracket $\left[x,y\right]=x\dashv y-y\vdash x$ turns every dialgebra into a Leibniz algebra, exactly as the commutator turns an associative algebra into a Lie algebra. Moreover, the universal enveloping dialgebra of a Leibniz algebra plays the same role as the universal enveloping associative algebra of a Lie algebra~\cite{Loday01}. Normal forms and Gr\"obner--Shirshov bases for dialgebras were studied in~\cite{BokutChenLiu}.

Kolesnikov~\cite{Kolesnikov08} (see also Pozhidaev~\cite{Pozhidaev08}) found a general procedure which associates with every variety $\mathcal{V}$ of algebras the class of $\mathcal{V}$-dialgebras: the two products $\dashv,\vdash$ are required to satisfy the so-called $0$-identities, and each defining identity of $\mathcal{V}$ is replaced by several identities obtained by ``emphasizing'' one of its variables. Applied to the variety of associative algebras, this procedure gives precisely Loday's dialgebras. Applied to the variety of (right) Leibniz algebras, it gives the class of Leibniz dialgebras (Definition~\ref{def:LD}), which is the subject of the present paper. Thus a Leibniz dialgebra is a Leibniz algebra $(D,\dashv)$ equipped with a second product $\vdash$ which is compatible with $\dashv$ in the sense of identities \eqref{LD2}--\eqref{LD5} below.

The algebraic classification of low-dimensional algebras of a given variety is a classical problem. Complex Leibniz algebras of dimension at most three are completely classified; we use the classification of three-dimensional complex Leibniz algebras due to Casas et al.~\cite{CILL} (Theorem~\ref{thm:list} below). For dialgebras, three-dimensional complex diassociative algebras were classified in~\cite{RRB10}, and four-dimensional nilpotent diassociative algebras in~\cite{BRR15}. Four-dimensional Leibniz dialgebras were considered in~\cite{NF4}, where all Leibniz dialgebras $(\mathcal{DL},\dashv,\vdash)$ such that $(\mathcal{DL},\dashv)$ is the four-dimensional null-filiform Leibniz algebra (see~\cite{AyupovOmirov01}) were described. Since for every Leibniz dialgebra $(\mathcal{DL},\dashv,\vdash)$ the algebra $(\mathcal{DL},\dashv)$ is a Leibniz algebra, the classification of three-dimensional Leibniz dialgebras reduces to the following problem: for each three-dimensional Leibniz algebra $(\mathcal{L},\dashv)$ of Theorem~\ref{thm:list}, find all products $\vdash$ on $\mathcal{L}$ such that $(\mathcal{L},\dashv,\vdash)$ is a Leibniz dialgebra, and classify them up to the action of the automorphism group $\operatorname{Aut}(\mathcal{L},\dashv).$

The first part of this classification was carried out in previous paper~\cite{KU}, where all three-dimensional complex Leibniz dialgebras $(\mathcal{DL},\dashv,\vdash)$ such that the Leibniz algebra $(\mathcal{DL},\dashv)$ has one-dimensional annihilator were described up to isomorphism. In this paper we solve this problem in the remaining non-Lie case $\dim\mathcal{L}^{\operatorname{ann}}=2,$ i.e. for the Leibniz algebras $\mathcal{L}_{8}\left(\alpha\right),$ $\mathcal{L}_9,$ $\mathcal{L}_{10},$ and $\mathcal{L}_{11}.$ The main results are Theorems~\ref{thm:L8}, \ref{thm:L9}, \ref{thm:L10} and~\ref{thm:L11}, summarized in Corollary~\ref{cor:all}: up to isomorphism there are $18$ families of such Leibniz dialgebras ($6$ over $\mathcal{L}_8(\alpha),$ $3$ over $\mathcal{L}_9,$ $6$ over $\mathcal{L}_{10}$ and $3$ over $\mathcal{L}_{11}$), nine of them depending on a complex parameter. Throughout the paper all vector spaces and algebras are finite-dimensional over the field $\mathbb{C}$ of complex numbers. In multiplication tables, products of basis vectors which are not written down are assumed to be zero.

\section{Preliminaries}

\begin{definition}\label{def:Leibniz}

A right Leibniz algebra is a vector space $\mathcal{L}$ with a bilinear bracket satisfying the Leibniz identity \begin{equation*} \left[\left[x,y\right],z\right]=\left[\left[x,z\right],y\right]+\left[x,\left[y,z\right]\right] \end{equation*} i.e. every right multiplication operator is a derivation. A Leibniz algebra with $[x,x]=0$ for all $x$ is a Lie algebra, and the ideal \begin{center} $\mathcal{L}^{\operatorname{ann}}=\mathrm{span}\{[x,x]\mid x\in\mathcal{L}\}$ \end{center} $($the annihilator of $\mathcal{L})$ is the smallest ideal of $\mathcal{L}$ such that $\mathcal{L}/\mathcal{L}^{\operatorname{ann}}$ is a Lie algebra.

\end{definition}

\begin{definition}\label{def:diassociative}

A dialgebra $($or diassociative algebra$)$ is a vector space $\mathcal{D}$ with two bilinear products $\dashv,\vdash\colon \mathcal{D}\times \mathcal{D}\to \mathcal{D}$ such 	that for all $x,y,z\in \mathcal{D}$ \begin{longtable}{ll}
$(x\dashv y)\dashv z=x\dashv(y\dashv z),$ & $(x\dashv y)\dashv z=x\dashv(y\vdash z),$ \\
$(x\vdash y)\dashv z=x\vdash(y\dashv z),$ & $(x\dashv y)\vdash z=(x\vdash y)\vdash z,$ \\
$(x\vdash y)\vdash z=x\vdash(y\vdash z).$ & \\
\end{longtable}

\end{definition}

\noindent For a dialgebra $(\mathcal{D},\dashv,\vdash)$ the bracket $\left[x,y\right]=x\dashv y-y\vdash x$ makes $\mathcal{D}$ a right Leibniz algebra (see~\cite{Loday01}).

Let $\mathcal{V}$ be a variety of algebras defined by multilinear identities. Following~\cite{Kolesnikov08}, a \textit{$\mathcal{V}$-dialgebra} is a vector space with two bilinear products $\dashv,\vdash$ satisfying \begin{enumerate}

\item the \textit{$0$-identities}: $\quad$ $x\dashv(y\dashv z)=x\dashv(y\vdash z),\qquad (x\dashv y)\vdash z=(x\vdash y)\vdash z;$

\item for every defining identity $f(x_1,\dots,x_n)=0$ of $\mathcal{V}$ and every $i=1,\dots,n,$ the identity $f^{(i)}(x_1,\dots,x_n)=0$ obtained from $f$ by emphasizing the variable $x_i:$ in each monomial of $f,$ every product $u\cdot v$ is replaced by $u\dashv v$ if $x_i$ occurs in $u$ and by $u\vdash v$ if $x_i$ occurs in $v$ (for products containing no $x_i$ the choice is irrelevant in view of the $0$-identities).

\end{enumerate} \noindent For the variety of associative algebras this gives exactly Definition~\ref{def:diassociative}. Applying the same procedure to the right Leibniz identity, emphasizing in turn $x,$ $y$ and $z,$ we arrive at the following definition, which is the central notion of this paper.

\begin{definition}\label{def:LD}

A Leibniz dialgebra is a vector space $\mathcal{DL}$ with two bilinear products $\dashv,\vdash\colon \mathcal{DL}\times \mathcal{DL}\to \mathcal{DL}$ such that for all $x,y,z\in \mathcal{DL}$ \begin{align}
(x\dashv y)\dashv z&=(x\dashv z)\dashv y+x\dashv(y\dashv z), \label{LD1}\tag{LD1} \\
(x\vdash y)\dashv z&=(x\dashv z)\vdash y+x\vdash(y\dashv z), \label{LD2}\tag{LD2} \\
(x\vdash y)\vdash z&=(x\vdash z)\dashv y+x\vdash(y\vdash z), \label{LD3}\tag{LD3} \\
x\dashv(y\dashv z)&=x\dashv(y\vdash z), \label{LD4}\tag{LD4} \\
(x\dashv y)\vdash z&=(x\vdash y)\vdash z. \label{LD5}\tag{LD5}
\end{align}

Identities \eqref{LD1}, \eqref{LD2} and \eqref{LD3} are the right Leibniz identity with the emphasized variable $x,$ $y$ and $z,$ respectively, and \eqref{LD4}, \eqref{LD5} are the $0$-identities.

\end{definition}

\begin{definition}\label{def:iso}

Two Leibniz dialgebras $(\mathcal{DL},\dashv,\vdash)$ and $(\mathcal{DL}^{\prime},\dashv^{\prime},\vdash^{\prime})$ are isomorphic if there is a linear bijection $\varphi\colon \mathcal{DL}\to \mathcal{DL}^{\prime}$ such that \begin{center} $\varphi(x\dashv y)=\varphi(x)\dashv^{\prime}\varphi(y),$ $\qquad$ $\varphi(x\vdash y)=\varphi(x)\vdash^{\prime}\varphi(y)$ $\qquad$ for all $x,y\in \mathcal{DL}.$ \end{center}

\end{definition}

\paragraph{Method of classification.} Let $(\mathcal{DL},\dashv,\vdash)$ be a three-dimensional Leibniz dialgebra such that $(\mathcal{DL},\dashv)$ is a non-Lie Leibniz algebra with $\dim\mathcal{L}^{\operatorname{ann}}=2.$ Then $(\mathcal{DL},\dashv)$ is isomorphic to one of the algebras of Theorem~\ref{thm:list}. Therefore we may assume that $(\mathcal{DL},\dashv)=(\mathcal{L},\dashv)$ is one of these algebras, with the basis $e_1,e_2,e_3$ and the multiplication table given there. Write \begin{center} $e_i\vdash e_j=\sum_{k=1}^{3}c_{ij}^{k}\,e_k, \qquad 1\le i,j\le 3,$ \end{center} with unknown structure constants $c_{ij}^k\in\mathbb{C}.$ Since the identities \eqref{LD1}--\eqref{LD5} are multilinear, it suffices to verify them for $x,y,z\in\{e_1,e_2,e_3\}.$ This gives a system of equations in the $c_{ij}^k$ whose solutions describe all Leibniz dialgebra structures with the given product $\dashv.$

An isomorphism of Leibniz dialgebras is, in particular, an isomorphism of the underlying Leibniz algebras. Hence dialgebras over non-isomorphic Leibniz algebras are non-isomorphic, and two structures $(\mathcal{L},\dashv,\vdash)$ and $(\mathcal{L},\dashv,\vdash^{\prime})$ are isomorphic if and only if there is an automorphism $\varphi\in\operatorname{Aut}(\mathcal{L},\dashv)$ such that $\varphi(x\vdash y)=\varphi(x)\vdash^{\prime}\varphi(y)$ for all $x,y\in \mathcal{L}.$ Thus the classification reduces to the description of the orbits of the group $\operatorname{Aut}(\mathcal{L},\dashv)$ on the set of solutions of the above system.

We use the classification of three-dimensional complex Leibniz algebras obtained in~\cite{CILL}; the algebras with two-dimensional annihilator are the following.

\begin{theorem}\label{thm:list}

Any three-dimensional complex non-Lie Leibniz algebra with $\dim\mathcal{L}^{\operatorname{ann}}=2$ is isomorphic to one of the following algebras:

\begin{longtable}{cll}
\hline
Algebra & Multiplication table & Automorphisms\\
\hline
\makecell{$\mathcal{L}_{8}\left(\alpha\right),$\\ $\alpha\in\mathbb{C}\setminus\{0\}$} & \makecell[l]{$\left[e_1, e_3\right]=\alpha e_1,$\\ $\left[e_2, e_3\right]=e_2.$} & see Remark~\ref{rem:L8}. \\
\hline
$\mathcal{L}_{9}$ & \makecell[l]{$\left[e_1, e_3\right]=e_1+e_2,$\\ $\left[e_2, e_3\right]=e_2.$} & \makecell[l]{$\varphi(e_1)=a_{11}e_1+a_{12}e_2,$\\ $\varphi(e_2)=a_{11}e_2,$\\ $\varphi(e_3)=e_3.$}\\
\hline
$\mathcal{L}_{10}$ & \makecell[l]{$\left[e_1, e_3\right]=e_2,$\\ $\left[e_3, e_3\right]=e_1.$} & \makecell[l]{$\varphi(e_1)=a_{33}^2e_1+a_{31}a_{33}e_2,$\\ $\varphi(e_2)=a_{33}^3e_2,$\\ $\varphi(e_3)=a_{31}e_1+a_{32}e_2+a_{33}e_3.$} \\
\hline
$\mathcal{L}_{11}$ & \makecell[l]{$\left[e_1, e_3\right]=e_2,$\\ $\left[e_2, e_3\right]=e_2,$\\ $\left[e_3, e_3\right]=e_1.$} & \makecell[l]{$\varphi(e_1)=e_1+(a_{22}-1)e_2,$\\ $\varphi(e_2)=a_{22}e_2,$\\ $\varphi(e_3)=(a_{22}-a_{32}-1)e_1+a_{32}e_2+e_3.$} \\
\hline
\end{longtable}

\noindent Here $a_{11}\neq0$ for $\mathcal{L}_{9},$ $a_{33}\neq0$ for $\mathcal{L}_{10}$ and $a_{22}\neq0$ for $\mathcal{L}_{11}.$

\end{theorem}

\begin{remark}\label{rem:L8} For $\alpha\neq0$ the algebras $\mathcal{L}_{8}\left(\alpha\right)$ and $\mathcal{L}_{8}\left(1/\alpha\right)$ are isomorphic. Indeed, in the new basis \begin{equation}\label{eq:swap} e_1^{\prime}=e_2, \quad e_2^{\prime}=e_1, \quad e_3^{\prime}=\alpha^{-1}e_3 \end{equation} of $\mathcal{L}_{8}\left(\alpha\right)$ we have $\left[e_1^{\prime},e_3^{\prime}\right]=\alpha^{-1}e_1^{\prime},$ $\left[e_2^{\prime},e_3'\right]=e_2^{\prime}.$ Conversely, let $\psi\colon\mathcal{L}_{8}\left(\alpha\right)\to\mathcal{L}_{8}\left(\beta\right)$ be an isomorphism. Then $\psi(\mathcal{L}^{\operatorname{ann}})=\mathcal{L}^{\operatorname{ann}}$ and $\psi(e_3)=c\,e_3+v$ with $c\neq0,$ $v\in\mathcal{L}^{\operatorname{ann}}.$ Since $\left[\mathcal{L}^{\operatorname{ann}},\mathcal{L}^{\operatorname{ann}}\right]=0,$ for $x\in\mathcal{L}^{\operatorname{ann}}$ we get $\psi\left(\left[x,e_3\right]\right)=\left[\psi(x),\psi(e_3)\right]=c\left[\psi(x),e_3\right],$ i.e. $\psi R_{e_3}=c\,R_{e_3}\psi$ on $\mathcal{L}^{\operatorname{ann}},$ where $R_{e_3}(x)=\left[x,e_3\right].$ Therefore the eigenvalues $\{c\beta,c\}$ of $cR_{e_3}$ in $\mathcal{L}_{8}\left(\beta\right)$ coincide with the eigenvalues $\{\alpha,1\}$ of $R_{e_3}$ in $\mathcal{L}_{8}\left(\alpha\right).$ Hence either $c=1,$ $\beta=\alpha,$ or $c=\alpha,$ $\beta=1/\alpha.$ Thus $\mathcal{L}_{8}\left(\alpha\right)\cong\mathcal{L}_{8}\left(\beta\right)$ iff $\beta\in\{\alpha,1/\alpha\}.$

The automorphism group of $\mathcal{L}_{8}\left(\alpha\right)$ is described as follows:

\begin{longtable}{lll}
\hline
$\alpha\neq\pm1$ & $:$ & $\varphi(e_1)=a_{11}e_1,$ $\varphi(e_2)=a_{22}e_2,$ $\varphi(e_3)=e_3,$
where $a_{11}a_{22}\neq0;$ \\
\hline
$\alpha=-1$ & $:$ & \makecell[l]{$\varphi(e_1)=a_{11}e_1,$ $\varphi(e_2)=a_{22}e_2,$ $\varphi(e_3)=e_3,$ where $a_{11}a_{22}\neq0,$\\ $\varphi(e_1)=a_{12}e_2,$ $\varphi(e_2)=a_{21}e_1,$ $\varphi(e_3)=-e_3,$ where $a_{12}a_{21}\neq0.$} \\
\hline
$\alpha=1$ & $:$ & \makecell[l]{$\varphi(e_1)=a_{11}e_1+a_{12}e_2,$ $\varphi(e_2)=a_{21}e_1+a_{22}e_2,$
$\varphi(e_3)=e_3,$\\ where $a_{11}a_{22}-a_{12}a_{21}\neq0;$} \\
\hline
\end{longtable}

\end{remark}

\section{Three-dimensional complex Leibniz dialgebras}

We choose the first part $\mathcal{A}_1=(\mathcal{DL},\dashv)$ of the Leibniz dialgebra from
Theorem~\ref{thm:list}, i.e. one of the algebras $\mathcal{L}_{8}\left(\alpha\right),$ $\mathcal{L}_{9},$ $\mathcal{L}_{10},$ $\mathcal{L}_{11}.$ The second part $\mathcal{A}_2=(\mathcal{DL},\vdash)$ is defined by the multiplication table \begin{equation}
\begin{aligned}
e_1 \vdash e_i &= \alpha_i e_1 + \beta_i e_2 + \gamma_i e_3, & i &\in \{1,2,3\}, \\
e_2 \vdash e_i &= \alpha_{i+3} e_1 + \beta_{i+3} e_2 + \gamma_{i+3} e_3, & i &\in \{1,2,3\}, \\
e_3 \vdash e_i &= \alpha_{i+6} e_1 + \beta_{i+6} e_2 + \gamma_{i+6} e_3, & i &\in \{1,2,3\}.
\end{aligned}
\label{eq:table}
\end{equation}

\noindent Imposing the Leibniz dialgebra axioms (Definition~\ref{def:LD}) we obtain constraints for the structure constants, and then we distinguish non-isomorphic algebras by the general change of basis (automorphisms of $\mathcal{A}_1$).

\begin{theorem}\label{thm:L8}

Any three-dimensional complex Leibniz dialgebra constructed from one of the algebras $\mathcal{L}_{8}\left(\alpha\right),$ $\alpha\in\mathbb{C}\setminus\{0\},$ is isomorphic to one of the following pairwise non-isomorphic algebras:

\begin{longtable}{lll}
\hline
$\mathcal{DL}_{8}^{1}\left(\alpha\right)$ & $:$ & $e_1\dashv e_3=\alpha e_1,$\ $e_2\dashv e_3=e_2.$ \\
\hline
$\mathcal{DL}_{8}^{2}\left(\alpha\right)$ & $:$ & $e_1\dashv e_3=\alpha e_1,$\ $e_2\dashv e_3=e_2,$\ $e_1\vdash e_3=\alpha e_1,$\ $e_2\vdash e_3=e_2.$ \\
\hline
$\mathcal{DL}_{8}^{3}\left(\alpha\right)$ & $:$ & $e_1\dashv e_3=\alpha e_1,$\ $e_2\dashv e_3=e_2,$\ $e_3\vdash e_1=-\alpha e_1,$\ $e_3\vdash e_2=-e_2.$ \\
\hline
\end{longtable}

where $\alpha\in\mathbb{C}\setminus\{0\}$ is considered up to $\alpha\mapsto1/\alpha$ $($one has $\mathcal{DL}_{8}^{i}\left(\alpha\right)\cong\mathcal{DL}_{8}^{i}\left(1/\alpha\right),$ $i=1,2,3)$

\begin{longtable}{lll}
\hline
$\mathcal{DL}_{8}^{4}\left(\alpha\right)$ & $:$ & $e_1\dashv e_3=\alpha e_1,$\ $e_2\dashv e_3=e_2,$\ $e_1\vdash e_3=\alpha e_1.$ \\
\hline
$\mathcal{DL}_{8}^{5}\left(\alpha\right)$ & $:$ & $e_1\dashv e_3=\alpha e_1,$\ $e_2\dashv e_3=e_2,$\ $e_3\vdash e_1=-\alpha e_1.$ \\
\hline
$\mathcal{DL}_{8}^{6}\left(\alpha\right)$ & $:$ & $e_1\dashv e_3=\alpha e_1,$\ $e_2\dashv e_3=e_2,$\ $e_1\vdash e_3=\alpha e_1,$\ $e_3\vdash e_2=-e_2.$ \\
\hline
\end{longtable}

where $\alpha\in\mathbb{C}\setminus\{0\}.$

\end{theorem}

\begin{proof}

Consider the algebra $\mathcal{A}_1=(\mathcal{DL},\dashv)$ with multiplication table \begin{center} $e_1\dashv e_3=\alpha e_1,$ $\quad$ $e_2\dashv e_3=e_2,$ $\quad$ $\alpha\in\mathbb{C}\setminus\{0\}.$ \end{center} The second part $\mathcal{A}_2=(\mathcal{DL},\vdash)$ is defined by the unknowns $\alpha_i,$ $\beta_i,$ $\gamma_i,$ $1\le i\le9,$ as in~\eqref{eq:table}. Imposing the Leibniz dialgebra axioms, we obtain \begin{longtable}{ll} $e_1\vdash e_3=\alpha_3e_1+\beta_3e_2,$ & $e_2\vdash e_3=\alpha_6e_1+\beta_6e_2,$ \\ $e_3\vdash e_1=\alpha_7e_1+\beta_7e_2,$ & $e_3\vdash e_2=\alpha_8e_1+\beta_8e_2,$ \end{longtable} \noindent and all the other products $e_i\vdash e_j$ are zero, with the following constraints: \begin{equation}\begin{array}{rrrr}
(\alpha-1)\beta_3=0, & (\alpha-1)\alpha_6=0, & (\alpha-1)\beta_7=0, & (\alpha-1)\alpha_8=0,\\[2pt]
\alpha_3^2+\alpha_6\beta_3=\alpha\alpha_3, & \beta_6^2+\alpha_6\beta_3=\beta_6, &
\beta_3(\alpha_3+\beta_6-\alpha)=0, & \alpha_6(\alpha_3+\beta_6-1)=0,\\[2pt]
\alpha_7^2+\alpha_8\beta_7=-\alpha\alpha_7, & \beta_8^2+\alpha_8\beta_7=-\beta_8, &
\beta_7(\alpha_7+\beta_8+1)=0, & \alpha_8(\alpha_7+\beta_8+\alpha)=0,\\[2pt]
\alpha_3\alpha_7+\alpha_6\beta_7=0, & \alpha_3\alpha_7+\beta_3\alpha_8=0, &
\beta_3\alpha_7+\beta_6\beta_7=0, & \alpha_3\beta_7+\beta_3\beta_8=0,\\[2pt]
\alpha_3\alpha_8+\alpha_6\beta_8=0, & \beta_3\alpha_8+\beta_6\beta_8=0, &
\alpha_6\alpha_7+\beta_6\alpha_8=0, & \alpha_6\beta_7+\beta_6\beta_8=0.
\end{array}\label{eq:L8cons}
\end{equation}

\noindent First we show that by a suitable change of basis one can always put $\beta_3=\alpha_6=\beta_7=\alpha_8=0.$ We have the following cases.

\begin{enumerate}

\item[Case 1.] Let $\alpha\neq1.$ Then the first row of \eqref{eq:L8cons} gives $\beta_3=\alpha_6=\beta_7=\alpha_8=0.$

\item[Case 2.] Let $\alpha=1.$ Let us consider the change of basis: \begin{equation*} \varphi(e_1)=a_{11}e_1+a_{12}e_2, \quad \varphi(e_2)=a_{21}e_1+a_{22}e_2, \quad \varphi(e_3)=e_3, \quad \Delta=a_{11}a_{22}-a_{12}a_{21}\neq0. \end{equation*} We write the new basis elements $\{\varphi(e_1),\varphi(e_2),\varphi(e_3)\}$ via the basis elements $\{e_1,e_2,e_3\}.$ By checking all the multiplications of the algebra in the new basis we obtain the relations between the parameters $\{\alpha_3^{\prime},\beta_3^{\prime},\alpha_6^{\prime},\beta_6^{\prime},\alpha_7^{\prime},\beta_7^{\prime},\alpha_8^{\prime},\beta_8^{\prime}\}$ and $\{\alpha_3,\beta_3,\alpha_6,\beta_6,\alpha_7,\beta_7,\alpha_8,\beta_8\}:$ from \begin{longtable}{ll} $\varphi(e_1)\vdash\varphi(e_3)=\alpha_3^{\prime}\varphi(e_1)+\beta_3^{\prime}\varphi(e_2),$ & $\varphi(e_2)\vdash\varphi(e_3)=\alpha_6^{\prime}\varphi(e_1)+\beta_6^{\prime}\varphi(e_2),$ \\ $\varphi(e_3)\vdash\varphi(e_1)=\alpha_7^{\prime}\varphi(e_1)+\beta_7^{\prime}\varphi(e_2),$ & $\varphi(e_3)\vdash\varphi(e_2)=\alpha_8^{\prime}\varphi(e_1)+\beta_8^{\prime}\varphi(e_2)$ \end{longtable} we get \begin{longtable}{lll} $\alpha_3^{\prime}$ & $=$ & $\dfrac{a_{22}(a_{11}\alpha_3+a_{12}\alpha_6)-a_{21}(a_{11}\beta_3+a_{12}\beta_6)}{\Delta},$ \\
$\beta_3^{\prime}$ & $=$ & $\dfrac{\beta_3a_{11}^2+(\beta_6-\alpha_3)a_{11}a_{12}-\alpha_6a_{12}^2}{\Delta},$ \\
$\alpha_6^{\prime}$ & $=$ & $-\dfrac{\beta_3a_{21}^2+(\beta_6-\alpha_3)a_{21}a_{22}-\alpha_6a_{22}^2}{\Delta},$ \\
$\beta_6^{\prime}$ & $=$ & $\dfrac{a_{11}(a_{21}\beta_3+a_{22}\beta_6)-a_{12}(a_{21}\alpha_3+a_{22}\alpha_6)}{\Delta},$ \\
$\alpha_7^{\prime}$ & $=$ & $\dfrac{a_{22}(a_{11}\alpha_7+a_{12}\alpha_8)-a_{21}(a_{11}\beta_7+a_{12}\beta_8)}{\Delta},$ \\
$\beta_7^{\prime}$ & $=$ & $\dfrac{\beta_7a_{11}^2+(\beta_8-\alpha_7)a_{11}a_{12}-\alpha_8a_{12}^2}{\Delta},$ \\
$\alpha_8^{\prime}$ & $=$ & $-\dfrac{\beta_7a_{21}^2+(\beta_8-\alpha_7)a_{21}a_{22}-\alpha_8a_{22}^2}{\Delta},$ \\
$\beta_8^{\prime}$ & $=$ & $\dfrac{a_{11}(a_{21}\beta_7+a_{22}\beta_8)-a_{12}(a_{21}\alpha_7+a_{22}\alpha_8)}{\Delta}.$
\end{longtable} In particular, \begin{equation} \alpha_3^{\prime}+\beta_6^{\prime}=\alpha_3+\beta_6, \quad \alpha_7^{\prime}+\beta_8^{\prime}=\alpha_7+\beta_8. \label{eq:traces} \end{equation} Note that the new parameters satisfy the same constraints \eqref{eq:L8cons} (with $\alpha=1$), since they are the structure constants of the same dialgebra in another basis.

\end{enumerate}

\begin{enumerate}

\item[2.1.] Let $\alpha_3+\beta_6\neq1$ and $\alpha_7+\beta_8\neq-1.$ Then the last two constraints of the second and of the third rows of \eqref{eq:L8cons} give $\beta_3=\alpha_6=\beta_7=\alpha_8=0.$

\item[2.2.] Let $\alpha_3+\beta_6=1.$ Then $\alpha_3^2+\alpha_6\beta_3=\alpha_3$ implies $\alpha_3\beta_6=\alpha_6\beta_3.$ We choose \begin{longtable}{lllll}
if $\beta_3\neq0:$ & $a_{11}=\alpha_3,$ & $a_{12}=\beta_3,$ & $a_{21}=\beta_6,$ & $a_{22}=-\beta_3;$ \\
if $\beta_3=0:$ & $a_{11}=1,$ & $a_{12}=0,$ & $a_{21}=\alpha_6,$ & $a_{22}=\beta_6-\alpha_3.$ \\
\end{longtable} \noindent In the first case $\Delta=-\beta_3(\alpha_3+\beta_6)=-\beta_3\neq0.$ In the second case $\alpha_3\beta_6=\alpha_6\beta_3=0,$ hence $(\beta_6-\alpha_3)^2=(\alpha_3+\beta_6)^2-4\alpha_3\beta_6=1$ and $\Delta=\beta_6-\alpha_3\neq0.$ In both cases a direct check using $\alpha_3\beta_6=\alpha_6\beta_3$ shows that \begin{equation*} \beta_3a_{11}^2+(\beta_6-\alpha_3)a_{11}a_{12}-\alpha_6a_{12}^2=0, \qquad \beta_3a_{21}^2+(\beta_6-\alpha_3)a_{21}a_{22}-\alpha_6a_{22}^2=0, \end{equation*} i.e. $\beta_3^{\prime}=\alpha_6^{\prime}=0.$ By \eqref{eq:traces} we have $\alpha_3^{\prime}+\beta_6^{\prime}=1.$ Now the constraints $\alpha_3^{\prime}\beta_7^{\prime}+\beta_3^{\prime}\beta_8^{\prime}=0$ and $\beta_3^{\prime}\alpha_7^{\prime}+\beta_6^{\prime}\beta_7^{\prime}=0$ give $\alpha_3^{\prime}\beta_7^{\prime}=\beta_6^{\prime}\beta_7^{\prime}=0,$ hence $\beta_7^{\prime}=(\alpha_3^{\prime}+\beta_6^{\prime})\beta_7^{\prime}=0.$ Similarly, $\alpha_3^{\prime}\alpha_8^{\prime}+\alpha_6^{\prime}\beta_8^{\prime}=0$ and $\alpha_6^{\prime}\alpha_7^{\prime}+\beta_6^{\prime}\alpha_8^{\prime}=0$ give $\alpha_8^{\prime}=(\alpha_3^{\prime}+\beta_6^{\prime})\alpha_8^{\prime}=0.$ Thus $\beta_3^{\prime}=\alpha_6^{\prime}=\beta_7^{\prime}=\alpha_8^{\prime}=0.$

\item[2.3.] Let $\alpha_3+\beta_6\neq1$ and $\alpha_7+\beta_8=-1.$ Then $\beta_3=\alpha_6=0,$ and \eqref{eq:L8cons} gives $\alpha_3^2=\alpha_3,$ $\beta_6^2=\beta_6,$ $\alpha_3\alpha_7=0,$ $\beta_6\beta_8=0.$ Since $\alpha_3,\beta_6\in\{0,1\}$ and $\alpha_3+\beta_6\neq1,$ we have either $\alpha_3=\beta_6=1$ or $\alpha_3=\beta_6=0.$ In the first case $\alpha_7=\beta_8=0,$ which contradicts $\alpha_7+\beta_8=-1.$ Hence $\alpha_3=\beta_6=0,$ i.e. $e_1\vdash e_3=e_2\vdash e_3=0,$ and then $\alpha_3^{\prime}=\beta_3^{\prime}=\alpha_6^{\prime}=\beta_6^{\prime}=0$ for every change of basis. Further, $\alpha_7^2+\alpha_8\beta_7=-\alpha_7$ and $\alpha_7+\beta_8=-1$ imply $\alpha_7\beta_8=\alpha_8\beta_7.$ We choose \begin{longtable}{lllll}
if $\beta_7\neq0:$ & $a_{11}=\alpha_7,$ & $a_{12}=\beta_7,$ & $a_{21}=\beta_8,$ & $a_{22}=-\beta_7;$ \\
if $\beta_7=0:$ & $a_{11}=1,$ & $a_{12}=0,$ & $a_{21}=\alpha_8,$ & $a_{22}=\beta_8-\alpha_7.$ \\
\end{longtable} \noindent In the first case $\Delta=-\beta_7(\alpha_7+\beta_8)=\beta_7\neq0.$ In the second case $(\beta_8-\alpha_7)^2=(\alpha_7+\beta_8)^2-4\alpha_7\beta_8=1$ and $\Delta=\beta_8-\alpha_7\neq0.$ As in the previous case we obtain $\beta_7^{\prime}=\alpha_8^{\prime}=0.$

Thus in all cases we may assume $\beta_3=\alpha_6=\beta_7=\alpha_8=0,$ and the constraints \eqref{eq:L8cons} take the form \begin{equation}\label{eq:L8diag}\begin{array}{lll}
\alpha_3(\alpha_3-\alpha)=0, & \beta_6(\beta_6-1)=0, & \alpha_7(\alpha_7+\alpha)=0, \\
\beta_8(\beta_8+1)=0, & \alpha_3\alpha_7=0, & \beta_6\beta_8=0. \\
\end{array}
\end{equation} Let us consider the change of basis \begin{equation*} \varphi(e_1)=a_{11}e_1, \quad \varphi(e_2)=a_{22}e_2, \quad \varphi(e_3)=e_3, \end{equation*} which is the general automorphism of $\mathcal{L}_{8}\left(\alpha\right)$ for $\alpha\neq\pm1.$ By checking all the multiplications of the algebra in the new basis we obtain the relations between the parameters $\{\alpha_3^{\prime},\beta_6^{\prime},\alpha_7^{\prime},\beta_8^{\prime}\}$ and $\{\alpha_3,\beta_6,\alpha_7,\beta_8\}:$ \begin{longtable}{lll} $\varphi(e_1)\vdash\varphi(e_3)=\alpha_3^{\prime}\varphi(e_1)$ & $\Rightarrow$ & $\alpha_3^{\prime}=\alpha_3,$ \\
$\varphi(e_2)\vdash\varphi(e_3)=\beta_6^{\prime}\varphi(e_2)$ & $\Rightarrow$ & $\beta_6^{\prime}=\beta_6,$ \\
$\varphi(e_3)\vdash\varphi(e_1)=\alpha_7^{\prime}\varphi(e_1)$ & $\Rightarrow$ & $\alpha_7^{\prime}=\alpha_7,$ \\
$\varphi(e_3)\vdash\varphi(e_2)=\beta_8^{\prime}\varphi(e_2)$ & $\Rightarrow$ & $\beta_8^{\prime}=\beta_8.$ \\
\end{longtable}The change of basis~\eqref{eq:swap}, which transforms $\mathcal{L}_{8}\left(\alpha\right)$ into $\mathcal{L}_{8}\left(1/\alpha\right)$ (see Remark~\ref{rem:L8}), gives \begin{equation} \alpha_3^{\prime}=\alpha^{-1}\beta_6, \qquad \beta_6^{\prime}=\alpha^{-1}\alpha_3, \qquad \alpha_7^{\prime}=\alpha^{-1}\beta_8, \qquad \beta_8^{\prime}=\alpha^{-1}\alpha_7. \label{eq:swaprule}
\end{equation} For $\alpha=-1$ the additional automorphism $\varphi(e_1)=a_{12}e_2,$ $\varphi(e_2)=a_{21}e_1,$ $\varphi(e_3)=-e_3$ is the composition of \eqref{eq:swap} with a change of the above form. By \eqref{eq:L8diag} we have $\alpha_3\in\{0,\alpha\},$ $\alpha_7\in\{0,-\alpha\},$ $\beta_6\in\{0,1\},$ $\beta_8\in\{0,-1\},$ $\alpha_3\alpha_7=0$ and $\beta_6\beta_8=0.$ Then we have the following cases. \begin{itemize}
\item[$(1)$] Let $\alpha_3=\alpha_7=0.$
\begin{itemize}
\item[$(1.1)$] If $\beta_6=\beta_8=0,$ then we obtain $\mathcal{DL}^{1}_{8}\left(\alpha\right).$
\item[$(1.2)$] If $\beta_6=1,$ $\beta_8=0,$ then we obtain $\mathcal{DL}^{4}_{8}\left(1/\alpha\right).$
\item[$(1.3)$] If $\beta_6=0,$ $\beta_8=-1,$ then we obtain $\mathcal{DL}^{5}_{8}\left(1/\alpha\right).$
\end{itemize}
\item[$(2)$] Let $\alpha_3=\alpha,$ $\alpha_7=0.$
\begin{itemize}
\item[$(2.1)$] If $\beta_6=\beta_8=0,$ then we obtain $\mathcal{DL}^{4}_{8}\left(\alpha\right).$
\item[$(2.2)$] If $\beta_6=1,$ $\beta_8=0,$ then we obtain $\mathcal{DL}^{2}_{8}\left(\alpha\right).$
\item[$(2.3)$] If $\beta_6=0,$ $\beta_8=-1,$ then we obtain $\mathcal{DL}^{6}_{8}\left(\alpha\right).$
\end{itemize}
\item[$(3)$] Let $\alpha_3=0,$ $\alpha_7=-\alpha.$
\begin{itemize}
\item[$(3.1)$] If $\beta_6=\beta_8=0,$ then we obtain $\mathcal{DL}^{5}_{8}\left(\alpha\right).$
\item[$(3.2)$] If $\beta_6=1,$ $\beta_8=0,$ then we obtain $\mathcal{DL}^{6}_{8}\left(1/\alpha\right).$
\item[$(3.3)$] If $\beta_6=0,$ $\beta_8=-1,$ then we obtain $\mathcal{DL}^{3}_{8}\left(\alpha\right).$
\end{itemize}
\end{itemize}

\end{enumerate}

Finally, let us show that the obtained algebras are pairwise non-isomorphic. By Remark~\ref{rem:L8}, an algebra over $\mathcal{L}_{8}\left(\alpha\right)$ can be isomorphic to an algebra over $\mathcal{L}_{8}\left(\beta\right)$ only if $\beta=\alpha$ or $\beta=1/\alpha,$ and every isomorphism of $\mathcal{L}_{8}\left(\alpha\right)$ onto $\mathcal{L}_{8}\left(1/\alpha\right)$ is the composition of \eqref{eq:swap} with an automorphism of $\mathcal{L}_{8}\left(\alpha\right).$ The parameters $(\alpha_3,\beta_6,\alpha_7,\beta_8)$ of the algebras $\mathcal{DL}^{1}_{8}\left(\alpha\right),\dots,\mathcal{DL}^{6}_{8}\left(\alpha\right)$ are $(0,0,0,0),$ \ $(\alpha,1,0,0),$ \ $(0,0,-\alpha,-1),$ \ $(\alpha,0,0,0),$ \ $(0,0,-\alpha,0),$ \ $(\alpha,0,0,-1),$ respectively.

Let $\alpha\neq1.$ Then the parameters $(\alpha_3,\beta_6,\alpha_7,\beta_8)$ do not change under the automorphisms of $\mathcal{L}_{8}\left(\alpha\right)$ with $\varphi(e_3)=e_3,$ and \eqref{eq:swap} transforms them by the rule \eqref{eq:swaprule}. Hence an algebra of the list over $\mathcal{L}_{8}\left(\alpha\right)$ is isomorphic to an algebra of the list over $\mathcal{L}_{8}\left(\beta\right)$ only if either $\beta=\alpha$ and their parameters coincide, or $\beta=1/\alpha$ and their parameters are related by \eqref{eq:swaprule} (for $\alpha=-1$ both possibilities occur with $\beta=\alpha$). The rule \eqref{eq:swaprule} transforms the parameters of $\mathcal{DL}^{i}_{8}\left(\alpha\right),$ $i=1,2,3,$ into the parameters $(0,0,0,0),$ $(\alpha^{-1},1,0,0),$ $(0,0,-\alpha^{-1},-1)$ of $\mathcal{DL}^{i}_{8}\left(1/\alpha\right),$ and transforms the parameters of $\mathcal{DL}^{4}_{8}\left(\alpha\right),$ $\mathcal{DL}^{5}_{8}\left(\alpha\right),$ $\mathcal{DL}^{6}_{8}\left(\alpha\right)$ into $(0,1,0,0),$ \ $(0,0,0,-1),$ \ $(0,1,-\alpha^{-1},0),$ which do not occur in the list over $\mathcal{L}_{8}\left(1/\alpha\right)$ (an algebra of the list with $\beta_6=1$ has $\alpha_3\neq0,$ and an algebra of the list with $\beta_8=-1$ has $\alpha_3\neq0$ or $\alpha_7\neq0$).

Let $\alpha=1.$ By Remark~\ref{rem:L8}, algebras over $\mathcal{L}_{8}\left(1\right)$ can be isomorphic only to algebras over $\mathcal{L}_{8}\left(1\right),$ and by \eqref{eq:traces} the numbers $\alpha_3+\beta_6$ and $\alpha_7+\beta_8$ are invariant under the general change of basis. For $\mathcal{DL}^{1}_{8}\left(1\right),\dots,\mathcal{DL}^{6}_{8}\left(1\right)$ the pairs $(\alpha_3+\beta_6,\alpha_7+\beta_8)$ are $(0,0),$ \ $(2,0),$ \ $(0,-2),$ \ $(1,0),$ \ $(0,-1),$ \ $(1,-1),$ respectively, and they are pairwise distinct.

\end{proof}

\begin{theorem}\label{thm:L9}

Any three-dimensional complex Leibniz dialgebra constructed from the algebra $\mathcal{L}_9$ is isomorphic to one of the following pairwise non-isomorphic algebras:

\begin{longtable}{lll}
\hline
$\mathcal{DL}_{9}^{1}$ & $:$ & $e_1\dashv e_3=e_1+e_2,$ \ $e_2\dashv e_3=e_2.$ \\
\hline
$\mathcal{DL}_{9}^{2}$ & $:$ & $e_1\dashv e_3=e_1+e_2,$ \ $e_2\dashv e_3=e_2,$ \ $e_1\vdash e_3=e_1+e_2,$ \ $e_2\vdash e_3=e_2.$ \\
\hline
$\mathcal{DL}_{9}^{3}$ & $:$ & $e_1\dashv e_3=e_1+e_2,$ \ $e_2\dashv e_3=e_2,$ \ $e_3\vdash e_1=-e_1-e_2,$ \ $e_3\vdash e_2=-e_2.$ \\
\hline
\end{longtable}

\end{theorem}

\begin{proof}

Consider the algebra $\mathcal{A}_1=(\mathcal{DL},\dashv)$ with multiplication table \begin{center} $e_1\dashv e_3=e_1+e_2,$ $\quad$ $e_2\dashv e_3=e_2.$ \end{center} The second part $\mathcal{A}_2=(\mathcal{DL},\vdash)$ is defined by the unknowns $\alpha_i,$ $\beta_i,$ $\gamma_i,$ $1\le i\le9,$ as in~\eqref{eq:table}. Imposing the Leibniz dialgebra axioms, we obtain \begin{longtable}{ll} $e_1\vdash e_3=\alpha_3e_1+\beta_3e_2,$ & $e_2\vdash e_3=\alpha_3e_2,$ \\ $e_3\vdash e_1=\alpha_7e_1+\beta_7e_2,$ & $e_3\vdash e_2=\alpha_7e_2,$ \end{longtable} \noindent and all the other products $e_i\vdash e_j$ are zero, with the following constraints: \begin{equation*}\begin{array}{rrr}
\alpha_3(\alpha_3-1)=0, & \beta_3(1-2\alpha_3)+\alpha_3=0, & \alpha_3\alpha_7=0,\\
\alpha_7(\alpha_7+1)=0, & \beta_7(1+2\alpha_7)+\alpha_7=0, & \alpha_7\beta_3+\alpha_3\beta_7=0.
\end{array}
\end{equation*}

\noindent Let us consider the change of basis \begin{equation*} \varphi(e_1)=a_{11}e_1+a_{12}e_2, \quad \varphi(e_2)=a_{11}e_2, \quad \varphi(e_3)=e_3, \quad a_{11}\neq0. \end{equation*} We write the new basis elements $\{\varphi(e_1),\varphi(e_2),\varphi(e_3)\}$ via the basis elements $\{e_1,e_2,e_3\}.$ By checking all the multiplications of the algebra in the new basis we obtain the relations between the parameters $\{\alpha_3^{\prime},\beta_3^{\prime},\alpha_7^{\prime},\beta_7^{\prime}\}$ and $\{\alpha_3,\beta_3,\alpha_7,\beta_7\}:$ \begin{longtable}{llll} $\varphi(e_1)\vdash\varphi(e_3)=\alpha_3^{\prime}\varphi(e_1)+\beta_3^{\prime}\varphi(e_2)$ & $\Rightarrow$ & $\alpha_3^{\prime}=\alpha_3,$ & $\beta_3^{\prime}=\beta_3,$ \\ $\varphi(e_3)\vdash\varphi(e_1)=\alpha_7^{\prime}\varphi(e_1)+\beta_7^{\prime}\varphi(e_2)$ & $\Rightarrow$ & $\alpha_7^{\prime}=\alpha_7,$ & $\beta_7^{\prime}=\beta_7.$ \\ \end{longtable} \noindent Then we have the following cases.

\begin{enumerate}

\item[Case 1.] Let $\alpha_3=\alpha_7=0.$ Then $\beta_3=\beta_7=0,$ and we obtain $\mathcal{DL}_{9}^{1}.$
\item[Case 2.] Let $\alpha_3=1.$ Then $\alpha_7=0,$ $\beta_3=1,$ $\beta_7=0,$ and we obtain $\mathcal{DL}_{9}^{2}.$
\item[Case 3.] Let $\alpha_7=-1.$ Then $\alpha_3=0,$ $\beta_3=0,$ $\beta_7=-1,$ and we obtain $\mathcal{DL}_{9}^{3}.$

\end{enumerate}

Since all the parameters are invariant under the change of basis, the algebras $\mathcal{DL}_{9}^{1},$ $\mathcal{DL}_{9}^{2},$ $\mathcal{DL}_{9}^{3}$ are pairwise non-isomorphic.

\end{proof}

\begin{theorem}\label{thm:L10}

Any three-dimensional complex Leibniz dialgebra constructed from the algebra $\mathcal{L}_{10}$ is isomorphic to one of the following pairwise non-isomorphic algebras:

\begin{longtable}{lll}
\hline
$\mathcal{DL}_{10}^{1}$ & $:$ & $e_1\dashv e_3=e_2,$ \ $e_3\dashv e_3=e_1.$ \\
\hline
$\mathcal{DL}_{10}^{2}$ & $:$ & $e_1\dashv e_3=e_2,$ \ $e_3\dashv e_3=e_1,$ \ $e_3\vdash e_3=e_2.$ \\
\hline
$\mathcal{DL}_{10}^{3}$ & $:$ & $e_1\dashv e_3=e_2,$ \ $e_3\dashv e_3=e_1,$ \ $e_1\vdash e_3=e_2,$ \ $e_3\vdash e_3=e_1.$ \\
\hline
$\mathcal{DL}_{10}^{4}$ & $:$ & $e_1\dashv e_3=e_2,$ \ $e_3\dashv e_3=e_1,$ \ $e_1\vdash e_3=e_2,$ \ $e_3\vdash e_3=e_1+e_2.$ \\
\hline
$\mathcal{DL}_{10}^{5}$ & $:$ & $e_1\dashv e_3=e_2,$ \ $e_3\dashv e_3=e_1,$ \ $e_3\vdash e_1=-e_2,$ \ $e_3\vdash e_3=-e_1.$ \\
\hline
$\mathcal{DL}_{10}^{6}$ & $:$ & $e_1\dashv e_3=e_2,$ \ $e_3\dashv e_3=e_1,$ \ $e_3\vdash e_1=-e_2,$ \ $e_3\vdash e_3=-e_1+e_2.$ \\
\hline
\end{longtable}

\end{theorem}

\begin{proof}

Consider the algebra $\mathcal{A}_1=(\mathcal{DL},\dashv)$ with multiplication table \begin{center} $e_1\dashv e_3=e_2,$ $\quad$ $e_3\dashv e_3=e_1.$ \end{center} The second part $\mathcal{A}_2=(\mathcal{DL},\vdash)$ is defined by the unknowns $\alpha_i,$ $\beta_i,$ $\gamma_i,$ $1\le i\le9,$ as in~\eqref{eq:table}. Imposing the Leibniz dialgebra axioms, we obtain \begin{longtable}{lll}
$e_1\vdash e_3=(\alpha_9-\beta_7)e_2,$ & $e_3\vdash e_1=\beta_7e_2,$ & $e_3\vdash e_3=\alpha_9e_1+\beta_9e_2,$ \end{longtable} \noindent and all the other products $e_i\vdash e_j$ are zero, with the following constraints: \begin{equation*} (\alpha_9-1)(\alpha_9-\beta_7)=0, \quad \beta_7(\alpha_9+1)=0, \end{equation*} i.e. $(\beta_7,\alpha_9)\in\{(0,0),(0,1),(-1,-1)\}$ and $\beta_9\in\mathbb{C}$ is arbitrary. Let us consider the change of the generators of basis \begin{equation*} \varphi(e_1)=a_{33}^2e_1+a_{31}a_{33}e_2, \quad \varphi(e_2)=a_{33}^3e_2, \quad \varphi(e_3)=a_{31}e_1+a_{32}e_2+a_{33}e_3, \quad a_{33}\neq0. \end{equation*} We write the new basis elements $\{\varphi(e_1),\varphi(e_2),\varphi(e_3)\}$ via the basis elements $\{e_1,e_2,e_3\}.$ By checking all the multiplications of the algebra in the new basis we obtain the relations between the parameters $\{\beta_7^{\prime},\alpha_9^{\prime},\beta_9^{\prime}\}$ and $\{\beta_7,\alpha_9,\beta_9\}:$ \begin{longtable}{lll} $\varphi(e_1)\vdash\varphi(e_3)=\left(\alpha_9^{\prime}-\beta_7^{\prime}\right)\varphi(e_2)$ & $\Rightarrow$ & $\alpha_9^{\prime}-\beta_7^{\prime}=\alpha_9-\beta_7,$ \\
$\varphi(e_3)\vdash\varphi(e_1)=\beta_7^{\prime}\varphi(e_2)$ & $\Rightarrow$ &  $\beta_7^{\prime}=\beta_7.$ \\
$\varphi(e_3)\vdash\varphi(e_3)=\alpha_9^{\prime}\varphi(e_1)+\beta_9^{\prime}\varphi(e_2)$ & $\Rightarrow$ & $\alpha_9^{\prime}=\alpha_9,$ $\quad$ $\beta_9^{\prime}=\frac{1}{a_{33}}\beta_9.$ \\
\end{longtable} \noindent Hence, if $\beta_9=0,$ then $\beta_9^{\prime}=0$ for every change of basis, and if $\beta_9\neq0,$ then by choosing $a_{33}=\beta_9$ we obtain $\beta_9^{\prime}=1.$ Then we have the following cases.

\begin{enumerate}

\item[Case 1.] Let $\beta_7=\alpha_9=0.$ If $\beta_{9}=0,$ then we obtain $\mathcal{DL}_{10}^{1};$ if $\beta_{9}\ne 0,$ then we obtain $\mathcal{DL}_{10}^{2}.$
\item[Case 2.] Let $\beta_7=0,$ $\alpha_9=1.$ If $\beta_{9}=0,$ then we obtain $\mathcal{DL}_{10}^{3};$ if $\beta_{9}\ne 0,$ then we obtain $\mathcal{DL}_{10}^{4}.$
\item[Case 3.] Let $\beta_7=\alpha_9=-1.$ If $\beta_{9}=0,$ then we obtain $\mathcal{DL}_{10}^{5};$ if $\beta_{9}\ne 0,$ then we obtain $\mathcal{DL}_{10}^{6}.$

\end{enumerate}

Since $\beta_7$, $\alpha_9$ are invariant and $\beta_9$ can only be multiplied by a non-zero coefficient, the algebras $\mathcal{DL}_{10}^{1},\dots,\mathcal{DL}_{10}^{6}$ are pairwise non-isomorphic.

\end{proof}

\begin{theorem}\label{thm:L11}

Any three-dimensional complex Leibniz dialgebra constructed from the algebra $\mathcal{L}_{11}$ is isomorphic to one of the following pairwise non-isomorphic algebras:

\begin{longtable}{lll}
\hline
$\mathcal{DL}_{11}^{1}\left(q\right)$ & $:$ & $e_1\dashv e_3=e_2,$ \ $e_2\dashv e_3=e_2,$ \ $e_3\dashv e_3=e_1,$ \ $e_3\vdash e_3=qe_1-qe_2.$ \\
\hline
$\mathcal{DL}_{11}^{2}\left(q\right)$ & $:$ & \makecell[l]{$e_1\dashv e_3=e_2,$ \ $e_2\dashv e_3=e_2,$ \ $e_3\dashv e_3=e_1,$ \\ $e_1\vdash e_3=e_2,$ \ $e_2\vdash e_3=e_2,$ \ $e_3\vdash e_3=qe_1+(1-q)e_2.$} \\
\hline
$\mathcal{DL}_{11}^{3}\left(q\right)$ & $:$ & \makecell[l]{$e_1\dashv e_3=e_2,$ \ $e_2\dashv e_3=e_2,$ \ $e_3\dashv e_3=e_1,$ \\ $e_3\vdash e_1=-e_2,$ \ $e_3\vdash e_2=-e_2,$ \ $e_3\vdash e_3=qe_1-(1+q)e_2.$} \\
\hline
\end{longtable} where $q\in\mathbb{C}.$

\end{theorem}

\begin{proof}

Consider the algebra $\mathcal{A}_1=(\mathcal{DL},\dashv)$ with multiplication table \begin{center} $e_1\dashv e_3=e_2,$ $\quad$ $e_2\dashv e_3=e_2,$ $\quad$ $e_3\dashv e_3=e_1.$ \end{center} The second part $\mathcal{A}_2=(\mathcal{DL},\vdash)$ is defined by the unknowns $\alpha_i,$ $\beta_i,$ $\gamma_i,$ $1\le i\le9,$ as in~\eqref{eq:table}. Imposing the Leibniz dialgebra axioms, we obtain \begin{longtable}{lll} $e_1\vdash e_3=\beta_6e_2,$ & $e_2\vdash e_3=\beta_6e_2,$ & \\ $e_3\vdash e_1=\beta_8e_2,$ & $e_3\vdash e_2=\beta_8e_2,$ & $e_3\vdash e_3=\alpha_9e_1+(\beta_6+\beta_8-\alpha_9)e_2,$ \end{longtable} \noindent and all the other products $e_i\vdash e_j$ are zero, with the following constraints: \begin{equation*} \beta_6(\beta_6-1)=0, \quad \beta_8(\beta_8+1)=0, \quad \beta_6\beta_8=0, \end{equation*} i.e. $(\beta_6,\beta_8)\in\{(0,0),(1,0),(0,-1)\}$ and $\alpha_9\in\mathbb{C}$ is arbitrary. Let us consider the change of the generators of basis \begin{equation*} \varphi(e_1)=e_1+(a_{22}-1)e_2, \quad \varphi(e_2)=a_{22}e_2, \quad \varphi(e_3)=(a_{22}-a_{32}-1)e_1+a_{32}e_2+e_3, \quad a_{22}\neq0. \end{equation*} We write the new basis elements $\{\varphi(e_1),\varphi(e_2),\varphi(e_3)\}$ via the basis elements $\{e_1,e_2,e_3\}.$ By checking all the multiplications of the algebra in the new basis we obtain the relations between the parameters $\{\beta_6^{\prime},\beta_8^{\prime},\alpha_9^{\prime}\}$ and $\{\beta_6,\beta_8,\alpha_9\}:$ \begin{longtable}{lll} $\varphi(e_1)\vdash\varphi(e_3)=\beta_6^{\prime}\varphi(e_2)$ & $\Rightarrow$ & $\beta_6^{\prime}=\beta_6,$ \\ $\varphi(e_3)\vdash\varphi(e_1)=\beta_8^{\prime}\varphi(e_2)$ & $\Rightarrow$ &  $\beta_8^{\prime}=\beta_8.$ \\ $\varphi(e_3)\vdash\varphi(e_3)=\alpha_9^{\prime}\varphi(e_1)+(\beta_6^{\prime}+\beta_8^{\prime}-\alpha_9^{\prime})\varphi(e_2)$ & $\Rightarrow$ & $\alpha_9^{\prime}=\alpha_9.$ \\ \end{longtable} \noindent Hence all the parameters are invariant, and the parameter $\alpha_9=q$ cannot be normalized. Then we have the following cases.

\begin{enumerate}

\item[Case 1.] Let $\beta_6=\beta_8=0.$ Then we obtain $\mathcal{DL}_{11}^{1}\left(q\right).$
\item[Case 2.] Let $\beta_6=1,$ $\beta_8=0.$ Then we obtain $\mathcal{DL}_{11}^{2}\left(q\right).$
\item[Case 3.] Let $\beta_6=0,$ $\beta_8=-1.$ Then we obtain $\mathcal{DL}_{11}^{3}\left(q\right).$

\end{enumerate}

Since $\beta_6,$ $\beta_8,$ $\alpha_9$ are invariant, all the algebras $\mathcal{DL}_{11}^{1}\left(q\right),$ $\mathcal{DL}_{11}^{2}\left(q\right),$ $\mathcal{DL}_{11}^{3}\left(q\right),$ $q\in\mathbb{C},$ are pairwise non-isomorphic.

\end{proof}

\begin{corollary}\label{cor:all}

Any three-dimensional complex Leibniz dialgebra $(\mathcal{DL},\dashv,\vdash)$ such that $(\mathcal{DL},\dashv)$ is a non-Lie Leibniz algebra with $\dim\mathcal{L}^{\operatorname{ann}}=2$ is isomorphic to exactly one of the algebras $\mathcal{DL}_{8}^{1}\left(\alpha\right),\dots,\mathcal{DL}_{8}^{6}\left(\alpha\right)$ of Theorem~\ref{thm:L8}, $\mathcal{DL}_{9}^{1},\mathcal{DL}_{9}^{2},\mathcal{DL}_{9}^{3}$ of Theorem~\ref{thm:L9}, $\mathcal{DL}_{10}^{1},\dots,\mathcal{DL}_{10}^{6}$ of Theorem~\ref{thm:L10} and $\mathcal{DL}_{11}^{1}\left(q\right),$ $\mathcal{DL}_{11}^{2}\left(q\right),$ $\mathcal{DL}_{11}^{3}\left(q\right)$ of Theorem~\ref{thm:L11} $($where in $\mathcal{DL}_{8}^{1}\left(\alpha\right)$, $\mathcal{DL}_{8}^{2}\left(\alpha\right)$, $\mathcal{DL}_{8}^{3}\left(\alpha\right)$ the parameter $\alpha$ is considered up to $\alpha\mapsto1/\alpha).$

\end{corollary}


\begin{thebibliography}{99}

\bibitem{AyupovOmirov01}
Sh.~A.~Ayupov, B.~A.~Omirov, On some classes of nilpotent Leibniz algebras, Sib. Math. J. 42 (2001), no.~1, 15--24.

\bibitem{AOR}
Sh.~Ayupov, B.~Omirov, I.~Rakhimov, Leibniz Algebras: Structure and Classification, CRC Press, Boca Raton, FL, 2020.

\bibitem{NF4}
M.~Azizov, Description of $4$-dimensional Leibniz dialgebras which are constructed by null-filiform Leibniz algebra, Bull. Inst. Math. 6 (2023), no.~4, 1--8.

\bibitem{BRR15}
W.~Basri, I.~S.~Rakhimov, I.~M.~Rikhsiboev, Four-dimensional nilpotent diassociative algebras, J. Gen. Lie Theory Appl. 9 (2015), no.~1, 1000218.

\bibitem{Bloh}
A.~Bloh, On a generalization of the concept of Lie algebra, Dokl. Akad. Nauk SSSR 165 (1965), 471--473.

\bibitem{BokutChenLiu}
L.~A.~Bokut, Y.~Chen, C.~Liu, Gr\"obner--Shirshov bases for dialgebras, Internat. J. Algebra Comput. 20 (2010), no.~3, 391--415.

\bibitem{CILL}
J.~M.~Casas, M.~A.~Insua, M.~Ladra, S.~Ladra, An algorithm for the classification of $3$-dimensional complex Leibniz algebras, Linear Algebra Appl. 436 (2012), no.~9, 3747--3756.

\bibitem{FKS}
R.~Fehlberg~J\'unior, I.~Kaygorodov, A.~Saydaliyev, The geometric classification of symmetric Leibniz algebras, Commun. Math. 33 (2025), no.~1, 10.

\bibitem{KhM}
A.~Khudoyberdiyev, Kh.~Muratova, Solvable Leibniz superalgebras whose nilradical has the characteristic sequence $(n-1,1\,\vert{}\,m)$ and nilindex $n+m$, Commun. Math. 32 (2024), no.~2, 27--54.

\bibitem{Kolesnikov08}
P.~S.~Kolesnikov, Varieties of dialgebras and conformal algebras, Sib. Math. J. 49 (2008), no.~2, 257--272.

\bibitem{KU}
T.~Kurbanbaev, N.~Uzakbaev, Classification of three-dimensional complex Leibniz dialgebras with one-dimensional annihilator, Uzbek Math. J. 70 (2026), no.~1, 116--125.

\bibitem{Loday93}
J.-L.~Loday, Une version non commutative des alg\`ebres de Lie: les alg\`ebres de Leibniz, Enseign. Math. (2) 39 (1993), no.~3--4, 269--293.

\bibitem{Loday95}
J.-L.~Loday, Alg\`ebres ayant deux op\'erations associatives (dig\`ebres), C. R. Acad. Sci. Paris S\'er. I Math. 321 (1995), no.~2, 141--146.

\bibitem{Loday01}
J.-L.~Loday, Dialgebras, in: Dialgebras and Related Operads, Lecture Notes in Math. 1763, Springer, Berlin, 2001, 7--66.

\bibitem{LodayPirashvili}
J.-L.~Loday, T.~Pirashvili, Universal enveloping algebras of Leibniz algebras and (co)homology, Math. Ann. 296 (1993), no.~1, 139--158.

\bibitem{Pozhidaev08}
A.~P.~Pozhidaev, Dialgebras and related triple systems, Sib. Math. J. 49 (2008), no.~4, 696--708.

\bibitem{RRB10}
I.~M.~Rikhsiboev, I.~S.~Rakhimov, W.~Basri, Classification of $3$-dimensional complex diassociative algebras, Malays. J. Math. Sci. 4 (2010), no.~2, 241--254.

\end{thebibliography}
\end{document}